%Paper: funct-an/9510001
%From: cortizo@ime.usp.br
%Date: Wed, 4 Oct 95 10:26:24-030

% ----------------------------------------------------------------
%   Plain TeX input file                Sep/1995
% ----------------------------------------------------------------

\magnification=\magstep1
\baselineskip=16 pt
\mathsurround=2 pt

% ----------------------------------------------------------------

\def\gr{\it}
\def\tit#1{\bigbreak\noindent{\bf #1}\medskip\nobreak}
\def\(#1){\smallskip\itemitem{(#1)}}
\def\iind(#1){\textindent{(#1)}}
\def\NU#1{\noindent\hbox to\parindent{{\bf#1}\hfil}\ignorespaces}
\def\PF{\noindent{\it Proof:\/} }
\def\FP{\hbox{\quad\vrule height6pt width4pt depth0pt}}

\def\0{\hbox{{\rm\O}}}
\def\sc{\subset}
\def\ct{\supset}
\def\X{\times}

\def\IFF{if and only if}
\def\NOT{\hbox{\rm not }}
\def\AND{{\rm\ and\ }}
\def\OR{{\rm\ or\ }}
\def\IMP{\Rightarrow}
\def\EQV{\Leftrightarrow}
\def\V{\,\forall\,}
\def\E{\,\exists\,}
\def\EU{\,\exists!\,}

\def\S#1{\Sigma(#1)}

\def\#{\equiv}
\def\<#1>{\langle#1\rangle}
\def\EXT#1{\,\overline{#1}\,}
\def\K#1{{\rm K}(#1)}

\def\R{{\bf R}}
\def\N{{\bf N}}
\def\Z{{\bf Z}}

\def\A{{\rm A}}
\def\B{{\rm B}}
\def\C{{\rm C}}
\def\D{{\rm D}}
\def\Ee{{\rm E}}
\def\En{\Ee^n}
\def\Ek{\Ee^k}

\def\G{{\rm G}}
\def\U{{\rm U}}
\def\Un{\U^n}
\def\Uk{\U^k}
\def\I{\rm I}
\def\Ve{\rm V}

\def\SA{\S\A}

\def\eq#1{{\rm eq}_{#1}}
\def\P{{\rm P}}
\def\Q{{\rm Q}}

\def\id#1{{\rm id}_{#1}}
\def\:{\colon}
\def\o{\circ}
\def\p#1#2{\pi^{#1}_{#2}}
\def\f{f}
\def\g{g}

\def\op{\oplus}
\def\od{\odot}

\def\RR{\EXT\R}
\def\ZZ{\EXT\Z}

\def\AA{\EXT\A}
\def\BB{\EXT\B}
\def\CC{\EXT\C}
\def\DD{\EXT\D}
\def\UU{\EXT\U}
\def\UUn{\UU^n}
\def\UUk{\UU^k}
\def\VV{\EXT\Ve}

\def\PP{\EXT\P}
\def\QQ{\EXT\Q}
\def\ff{\EXT\f}
\def\gg{\EXT\g}

\def\aa{\EXT a}
\def\bb{\EXT b}
\def\xx{\EXT x}

\def\KA{\K\A}
\def\KB{\K\B}
\def\KUn{\K\Un}

\def\a{\alpha}
\def\b{\beta}
\def\ga{\gamma}
\def\x{\xi}
\def\y{\upsilon}
\def\z{\zeta}

\def\8{\infty}

% ----------------------------------------------------------------

\hrule height0pt
\vfil
{\bf VIRTUAL EXTENSIONS}
\bigskip\bigskip\bigskip
{\bf S\'ergio F.~Cortizo}
\bigskip\bigskip
Instituto de Matem\'atica e Estat\'\i stica, Universidade de S\~ao
Paulo

Cidade Universit\'aria, Rua do Mat\~ao, 1010

05508-900, S\~ao Paulo, SP, Brasil
\smallskip
cortizo@ime.usp.br
\bigskip\bigskip\bigskip
{\bf Abstract}
\bigskip
\item{}A process of extending sets which can be used as foundation for an
alternative or\-gan\-ization for Differential and Integral Calculus is
presented.
\smallskip
PACS \ 02.90.+p
\vfil\eject

\tit{I. Introduction}

Our goal is to present an extension process that can be applied to any set.
It is a rel\-at\-ively simple construction, which can be developed strictly
inside the limits of Elementary Set Theory.

When applied to the ordered field $\R$ of real numbers, this process
introduces infinitesimal and infinite quantities, which can be used in an
alternative organization of Differential and Integral Calculus. This
application was our original motivation, and will be presented in a
subsequent work.

In Secs.~II, III and~IV we define concepts occurring in the statement of the
fundamental result of this work: the theorem presented in Sec.~V. In the
remainder sections we discuss that result.

\tit{II. Extension of Sets and Subsets}

Let $\A$ be any set. We will denote the set of all infinite sequences of
elements of~$\A$ by~$\SA$:
$$
   \SA = \left\{ (a_1,a_2,a_3,\ldots)\ |\ a_i\in\A,\ i\in\N \right\},
$$
where $\N=\{1,2,3,\ldots\}$ is the set of natural numbers.

We will introduce now a relation on the set~$\SA$: we will say that two
sequences $(a_i)=(a_1,a_2,\ldots)\in\SA$ and $(b_i)=(b_1,b_2,\ldots)\in\SA$
$(a_i)$ and~$(b_i)$ {\gr end equals\/} when there exists $n\in\N$ such that
$i>n$ implies $a_i=b_i$. This is an equivalence relation on~$\SA$, which we
will represent by `$\#$'. The equivalence class of~$(a_i)\in\SA$ will be
denoted by~$\<a_i>$, so $\<a_i>=\<b_i>$ \IFF\ $(a_i)\#(b_i)$.

The quotient $\AA=\SA/\#$ will be called {\gr virtual extension\/} of
set~$\A$, or simply {\gr extension of~$\A$}. In other words, the members
of~$\AA$ are the equivalence classes modulo~$\#$:
$$
   \AA = \left\{ \<a_i>\ |\ (a_i)\in\SA \right\}.
$$

Let now $\B\sc\A$ be any subset of~$\A$. We will say that a sequence
$(a_i)\in\SA$ {\gr ends in~$\B$} when $a_i\in\B$ after a certain value for
the index, i.e., when there exists $n\in\N$ such that $a_i\in\B$ for all
$i>n$. It is clear that if a sequence ends in~$\B$ then any other equivalent
sequence (by~$\#$) will also end in~$\B$. So, we can define the subset
$\BB\sc\AA$ of all classes $\<a_i>\in\AA$ whose representatives sequences end
in~$\B$:
$$
   \BB = \left\{ \<a_i>\in\AA\ |\ (a_i) \hbox{ ends in } \B\sc\A \right\}.
$$

Example: we will call {\gr virtual numbers}, or just {\gr virtuals}, the
members of the extension $\RR$ of the real numbers set~$\R$. Since $\Z$ is a
subset of~$\R$, we have the virtual extension $\ZZ\sc\RR$, whose elements
will be called {\gr virtual integers}. The members of~$\ZZ$ are represented
by sequences that assume, after a certain value of the index, only integer
values. An example of virtual integer is the class of the sequence
$(1,2,3,\ldots)\in\S\R$, which will be denoted simply by~$\8\in\ZZ$.

For any $a\in\A$, we will represent by $\aa\in\AA$ the equivalence class of
the sequence $(a,a,a,\ldots)\in\SA$ constant at~$a\in\A$. Besides, for any
subset $\B\sc\A$, we will denote by $\KB\sc\AA$ the class of all constant
sequences in~$\B$:
$$
   \KB = \left\{ \bb=\<b,b,b,\ldots>\in\AA\ |\ b\in\B \right\}.
$$

It is easy to see that:
\medskip
\NU{II.1}{\gr For any subset $\B\sc\A$, we have $\KB=\KA\cap\BB$.}
\medskip
\NU{II.2}{\gr For any subset $\B\sc\A$, we have:

\(i) $\BB=\AA$ \IFF\ $\B=\A${\rm;}

\(ii) $\BB=\0$ \IFF\ $\B=\0${\rm;}

\(iii) $\BB$ is unitary \IFF\ $\B$ is unitary.}
\medskip
\NU{II.3}{\gr If $\B$ and~$\C$ are two subsets of~$\A$ then:

\(i) $\BB\sc\CC$ \IFF\ $\B\sc\C${\rm;}

\(ii) $\BB=\CC$ \IFF\ $\B=\C$.}
\medskip

\tit{III. Relations and Functions}

The objective of this section is to establish terminology and notation. Many
definitions below are universal, but not all of them.

Let $\Ee$ be any set. We will identify a relation between $n$ variables
$x_i\in\Ee$ ($i=1,\ldots,n$) with the class of $n$-tuples
$(x_1,\ldots,x_n)\in\En$ which satisfy that relation, i.e., we are
considering an {\gr $n$-ary relation\/} on~$\Ee$ as a subset $\P\sc\En$ of
the Cartesian product $\En$ of $n$ copies of~$\Ee$. For example, the subset
$\eq\Ee\sc\Ee\X\Ee$ below is the {\gr equality relation\/} on~$\Ee$:
$$
    \eq\Ee = \left\{ (x,y)\in\Ee\X\Ee\ |\ x=y \right\}.
$$

We will often write
$$
    \P(x_1,x_2,\ldots,x_n)
$$
instead of
$$
    (x_1,x_2,\ldots,x_n)\in\P
$$
when a set $\P\sc\En$ is being interpreted as an $n$-ary relation. We will
also abbreviate the $n$-tuple `$(x_1,x_2,\ldots,x_n)$' to `$x$', writing
simply `$\P(x)$'.

For every $n$-ary relation $\P\sc\En$, we will denote its negation by
$(\NOT\P)\sc\En$, i.e.:

\(i) $(\NOT\P) = \left\{ x\in\En\ |\ x\notin\P \right\}$.\smallskip

\noindent Moreover, for every pair of $n$-ary relations $\P\sc\En$
and~$\Q\sc\En$ (the same $n\in\N$), we will use the following notation:

\(ii) $(\P \AND \Q) = \P \cap \Q$;

\(iii) $(\P \OR \Q) = \P \cup \Q$;

\(iv) $(\P \IMP \Q) = [(\NOT\P) \OR \Q]$;

\(v) $(\P \EQV \Q) = [(\P \IMP \Q) \AND (\Q \IMP \P)]$.\smallskip

Given a $(k+n)$-ary relation $\P\sc\Ee^{k+n}$, with $n\ge1$, we can fix its
$k$ first entries, leaving the remaining ones free, and thus construct an
$n$-ary relation between elements of the same set $\Ee$. If $a\in\Ek$, then
we define:
$$
    \P a = \left\{ x\in\En\ |\ \P(a,x) \right\}.
$$
It is clear that $\P a\sc\En$ is an $n$-ary relation on the same set~$\Ee$.
With this definition, the condition `$\P a(x)$' is equivalent to `$\P(a,x)$'.

We will use the symbol `$\V$' as abbreviation of `for every' (universal logic
quantifier), the symbol `$\E$' will mean `there exists' (existential
quantifier), and `$\EU$' will be an abbreviation for `there exists one and
only one'. For example, the statement:
$$
    \hbox{there exists $k\in\N$ such that, for all $i>k$, we have $\P(x_i)$}
$$
will be shortened to:
$$
    \E k\in\N,\V i>k,\ \P(x_i).
$$
Furthermore, those three symbols will be used to indicate the relations
constructed quant\-ifying the first entries of a given relation, according to
the following definitions.

If $\P\sc\Ee^{k+n}$ is a $(k+n)$-ary relation ($n\ge1$), and $\D\sc\Ek$ a
subset of~$\Ek$, then we define:

\(i) $(\V\D,\P) = \left\{ x\in\En\ |\  \V y\in\D,\ \P(y,x) \right\};$

\(ii) $(\E\D,\P) = \left\{ x\in\En\ |\  \E y\in\D,\ \P(y,x) \right\};$

\(iii) $(\EU\D,\P) = \left\{ x\in\En\ |\ \EU y\in\D,\ \P(y,x) \right\}.$
\smallskip

\noindent Obviously, $(\V\D,\P)\sc\En$ and $(\E\D,\P)\sc\En$, as well
as~$(\EU\D,\P)\sc\En$, are $n$-ary relations on the same set~$\Ee$, for which
we have:

\(i) $(\V\D,\P)(x)$ \IFF\ $[\V y\in\D,\ \P(y,x)$];

\(ii) $(\E\D,\P)(x)$ \IFF\ $[\E y\in\D,\ \P(y,x)$];

\(iii) $(\EU\D,\P)(x)$ \IFF\ $[\EU y\in\D,\ \P(y,x)$].\smallskip

Let $\P\sc\Ee^{n+k}$ be any $(n+k)$-ary relation. We will say that $\P$ is a
{\gr functional relation\/} from $\En$ to $\Ek$ if, for every $x\in\En$,
there exists at the most one $y\in\Ek$ such that $\P(x,y)$. Besides, given
$\D\sc\En$ and $\C\sc\Ek$, we will say that a functional relation
$\P\sc\Ee^{n+k}$ is a {\gr function from $\D$ into~$\C$} when, for every
$x\in\D$, there exists $y\in\C$ such that $\P(x,y)$. As usual, we will
indicate that $\f$ is a function from $\D$ into $\C$ by writing:
$$
    \f\:\D\to\C.
$$
We will say that two functions $\f\:\D_\f\to\C_\f$ and $\g\:\D_\g\to\C_\g$
form a {\gr chain\/} when $\C_\f\sc\D_\g$; in this case we define the {\gr
composite function\/} $(\g\o\f)\:\D_\f\to\C_\g$ by $(\g\o\f)(x)=\g[\f(x)]$,
for every $x\in\D_\f$. We will represent the {\gr identity function\/} on any
set $\D$ by $\id\D$, i.e., $(\id\D)(x)=x$, for all $x\in\D$.

The composition operation can also be done between a relation and a function:
if $\P\sc\Ek$ and $\f\:\D\to\C$, with $\D\sc\En$ and $\C\sc\Ek$, then
$(\P\o\f)\sc\En$ is the $n$-ary relation defined by:
$$
   (\P\o\f)(x) \qquad\hbox{\IFF}\qquad x\in\D \hbox{\quad and\quad}\P[\f(x)].
$$

If $\f_i\:\D\to\C$ are functions defined from the same domain into the same
set ($i=1,\ldots,n$), then we will denote the ``aggregation'' of the $n$
functions~$\f_i$ by
$$
    (\f_1,\ldots,\f_n)\:\D\to\C^n,
$$
which means:
$$
    (\f_1,\ldots,\f_n)(x) = (\f_1(x),\ldots,\f_n(x)).
$$

In addition, for every subset $\D\sc\En$, we will represent the function that
simply selects its $i$-th entry by
$$
    \p i\D\:\D\to\Ee.
$$
That means:
$$
    \p i\D(x_1,\ldots,x_n) = x_i \qquad (i=1,\ldots,n).
$$

\tit{IV. Extension of Relations and Functions}

Let $\U$ be a fixed non-empty set. We will now apply the virtual extension
process defined in Sec.~II simultaneously to $\U$ and to the product $\Un$,
thus obtaining the sets $\UU$ and~$\EXT{\Un}$ respectively. In principle, the
extension $\EXT{\Un}$ of the product is different from the product $\UUn$ of
the extensions, but in this work we will adopt a practically universal
identification in mathematics: a sequence of members of the product $\Un$ is
the same as the corresponding $n$-tuple of sequences in~$\U$. In other words,
we will consider $\S\Un=[\S\U]^n$. Thus, it is easily seen that two sequences
in the product will end equal (Sec.~II) \IFF\ its $n$ component sequences end
equal, for we are working with finite products only. We then have
$\EXT{\Un}=\UUn$, for every $n\in\N$.

Let now $\P\sc\Un$ be a generic relation between $n$ variables $x_i\in\U$
($i=1,\ldots,n$). According to the extension process of subsets defined in
Sec.~II, the virtual extension of the set $\P\sc\Un$ is a subset of
$\PP\sc\EXT{\Un}$. With the identification above, this extended subset
$\PP\sc\UUn=\EXT{\Un}$ defines a new relation between $n$ variables $\x_i$
which range over the extended set~$\UU$ ($i=1,\ldots,n$). We will call $\PP$
the {\gr virtual extension of relation~$\P$}, or simply the {\gr extension
of~$\P$}.

For instance, let `$<$' be the ordering relation between real numbers.
According to the construction above, we have defined a binary relation
`$\EXT<$' between virtual numbers. Thinking about the members of~$\RR$ as
classes of sequences, we have $\<a_i>\EXT<\<b_i>$ \IFF\ there exists $n\in\N$
such that $a_i<b_i$ for every $i>n$. Then, $\xx\EXT<\8$ for every $x\in\R$.

It is important to note that the virtual extension of $\P\sc\Un$ does not
depend on the interpretation of $\P$ as a subset or a relation. In other
words, if
$$
   \A = \left\{ x\in\Un\ |\ \P(x) \right\}
$$
then
$$
   \AA = \left\{ \x\in\UUn\ |\ \PP(\x) \right\}.
$$

For example, if $\R^*_+\sc\R$ is the set of real non-negative numbers:
$$
   \R^*_+ = \left\{ x\in\R\ |\ x>0 \right\},
$$
then
$$
   \EXT{\R^*_+} = \left\{ \x\in\RR\ |\ \x\EXT>\EXT0 \right\}.
$$

It is easy to see that $\P\sc\U^{n+k}$ is a functional relation \IFF\ its
virtual extension $\PP\sc\UU^{n+k}$ is a functional relation with respect to
the same entries. Thus, if $\f\sc\U^{n+k}$ is a function from $\D\sc\Un$
into~$\C\sc\Uk$:
$$
    \f\:\D\to\C,
$$
then its virtual extension $\ff\sc\UU^{n+k}$ is a function from the extension
$\DD\sc\UUn$ of the domain of~$\f$ into the extension $\CC\sc\UUk$ of the
counterdomain of~$\f$, i.e.:
$$
    \ff\:\DD\to\CC.
$$

If $\<x_i>\in\DD\sc\UUn$ then the sequence $(x_i)\in\S\Un$ ends
in~$\D\sc\Un$. Hence, we can evaluate $f(x_i)$ for every~$i$ greater than
certain $k\in\N$. It is not difficult to verify that
$\ff(\<x_i>)\in\CC\sc\UUk$ is the class of sequences $(y_i)\in\S\Uk$ for
which there exists $l>k$ such that $y_i=f(x_i)$ for every $i>l$.

Examples: for $\U=\R$, we have the operations $+\:\R^2\to\R$ and
$\X\:\R^2\to\R$, whose extensions are:
$$
    \EXT+\:\RR^2\to\RR \qquad\hbox{and}\qquad \EXT\X\:\RR^2\to\RR.
$$
If $\a\in\RR$ and $\b\in\RR$ are two generic virtual numbers, we can
calculate the sum $\a\EXT+\b\in\RR$ and the product $\a\EXT\X\b\in\RR$. So,
the symbols `$\8\EXT+\EXT1$' and `$\8\EXT\X(\8\EXT+\EXT1)$' represent two
well defined virtual numbers. Moreover, we have the function
$\ln\:\R_+^*\to\R$, whose extension is:
$$
    \EXT\ln\:\EXT{\R_+^*}\to\RR.
$$
Since $\8\EXT>\EXT0$, there exists $\EXT\ln\8\in\RR$.

We can apply the logical operations defined in Sec.~III both to relations
on~$\U$ and to its extensions on~$\UU$. Thus, $\EXT{\P\OR\Q}\sc\UUn$ is the
virtual extension of relation $(\P\OR\Q)\sc\Un$, while $(\PP\OR\QQ)\sc\UUn$
is the logical disjunction of the extensions $\PP\sc\UUn$ and~$\QQ\sc\UUn$.

The following section is a theorem which states the ``commutation rules''
between the virtual extension process and the operations on relations and
functions defined in the last section.

\tit{V. The Virtual Extension Theorem}

For any set $\U$, and for every $n$-ary relation $\P\sc\Un$, the condition
`$\P(x)$' is equi\-val\-ent to~`$\PP(\xx)$', i.e., the statement `$\P(x)$' is
true for certain $x\in\Un$ \IFF\ `$\PP(\xx)$' is true for
$\xx\in\KUn\sc\UUn$.

We also have:

\(i) $\EXT{\eq\U} = \eq\UU$;

\(ii) $\EXT{\NOT\P}\sc(\NOT\PP)$;

\(iii) $\EXT{\P\AND\Q} = (\PP\AND\QQ)$;

\(iv) $\EXT{\P\OR\Q}\ct(\PP\OR\QQ)$;

\(v) $\EXT{\P\IMP\Q}\sc(\PP\IMP\QQ)$;

\(vi) $\EXT{\P\EQV\Q}\sc(\PP\EQV\QQ)$;

\(vii) $\EXT{\P a} = \PP\aa$;

\(viii) $\EXT{(\V\D,\P)} = (\V\DD,\PP)$;

\(ix) $\EXT{(\E\D,\P)} = (\E\DD,\PP)$;

\(x) $\EXT{(\EU\D,\P)} = (\EU\DD,\PP)$;

\(xi) $\EXT{\P\o\f} = \PP\o\ff$;

\(xii) $\EXT{\g\o\f} = \gg\o\ff$ \qquad and\qquad $\EXT{\id\D} = \id\DD$;

\(xiii) $\EXT{(\f_1,\ldots,\f_n)} = (\EXT{\f_1},\ldots,\EXT{\f_n})$
        \qquad and\qquad $\EXT{\p i\D} = \p i\DD$.\smallskip

Furthermore, for any $\D\sc\Un$, we have:

(a) The statements `$\V x\in\D,\ \P(x)$' and `$\V\x\in\DD,\ \PP(\x)$' are
logically equivalent, i.e., the first is true \IFF\ the second is true;

(b) The statements `$\E x\in\D,\ \P(x)$' and `$\E\x\in\DD,\ \PP(\x)$' are
logically equivalent;

(c) The statements `$\EU x\in\D,\ \P(x)$' and `$\EU\x\in\DD,\ \PP(\x)$' are
logically equivalent.

The remainder of this section will be dedicated to the proof of the~VET
(Virtual Extension Theorem):

The first assertion is exactly the proposition~II.1 applied to subset
$\P\sc\Un$:
$$
    \K\P = \K\Un\cap\PP.
$$

For the first series of items, we have:

\iind(i) For any classes $\<x_i>\in\UU$ and~$\<y_i>\in\UU$, the following
statements are equivalent:
$$\eqalign{
  &\EXT{\eq\U}(\<x_i>,\<y_i>)\cr
  &\E k,\V i>k,\ \eq\U(x_i,y_i)\cr
  &\E k,\V i>k,\ x_i=y_i\cr
  &\<x_i>=\<y_i>\cr
  &\eq\UU(\<x_i>,\<y_i>).\cr
}$$

\iind(ii) If $\<x_i>\in\UUn$ then the condition `$\EXT{\NOT\P}(\<x_i>)$' is
equivalent to:
$$
    \E k,\V i>k,\ \NOT\P(x_i),
$$
which is sufficient for the validity of
$$
    \V k,\E i>k,\ \NOT\P(x_i),
$$
which, in turn, is equivalent to~`$\NOT\PP(\<x_i>)$'.

\iind(iii) It is enough to note that the following statements are
equivalent:
$$\eqalign{
  &(\PP\AND\QQ)(\<x_i>)\cr
  &\PP(\<x_i>)\AND\QQ(\<x_i>)\cr
  &[\E k_1,\V i>k_1,\ \P(x_i)]\AND[\E k_2,\V i>k_2,\ \Q(x_i)]\cr
  &\E k,\V i>k,\ [\P(x_i)\AND\Q(x_i)]\cr
  &\E k,\V i>k,\ (\P\AND\Q)(x_i)\cr
  &\EXT{\P\AND\Q}(\<x_i>).\cr
}$$

\iind(iv) First, we have the equivalences:
$$\eqalign{
  &(\PP\OR\QQ)(\<x_i>)\cr
  &\PP(\<x_i>)\OR\QQ(\<x_i>)\cr
  &[\E k_1,\V i>k_1,\ \P(x_i)]\OR[\E k_2,\V i>k_2,\ \Q(x_i)].\cr
}$$
This last statement implies the first below, which is equivalent to the
following ones:
$$\eqalign{
  &\E k,\V i>k,\ \P(x_i)\OR\Q(x_i)\cr
  &\E k,\V i>k,\ (\P\OR\Q)(x_i)\cr
  &\EXT{\P\OR\Q}(\<x_i>).\cr
}$$

\iind(v) We will prove that $[\NOT(\PP\IMP\QQ)]\sc[\NOT\EXT{\P\IMP\Q}]$. The
following statements are equivalent:
$$\eqalign{
  &[\NOT(\PP\IMP\QQ)](\<x_i>)\cr
  &[\PP\AND(\NOT\QQ)](\<x_i>)\cr
  &\PP(\<x_i>)\AND(\NOT\QQ)(\<x_i>)\cr
  &[\E k_1,\V i>k_1,\ \P(x_i)]\AND[\V k_2,\E i>k_2,\ \NOT\Q(x_i)].\cr
}$$
{}From this last assertion we conclude the first one below, which is
equivalent to the following:
$$\eqalign{
  &\V k,\E i>k,\ [\P(x_i)\AND\NOT\Q(x_i)]\cr
  &\V k,\E i>k,\ \NOT(\P\IMP\Q)(x_i)\cr
  &\NOT\EXT{\P\IMP\Q}(\<x_i>).\cr
}$$

\iind(vi) Applying items (iii) and~(v) above, we have:
$$\eqalign{
  \EXT{\P\EQV\Q} &= \EXT{(\P\IMP\Q)\AND(\Q\IMP\P)}\cr
                 &= [(\EXT{\P\IMP\Q})\AND(\EXT{\Q\IMP\P})]\cr
                 &\sc [(\PP\IMP\QQ)\AND(\QQ\IMP\PP)]\cr
                 &= (\PP\EQV\QQ).\cr
}$$

\iind(vii) It is enough to note that, if $a\in\U$ and~$\<x_i>\in\UUn$, then
the next following statements are equivalent:
$$\eqalign{
  &\EXT{\P a}(\<x_i>)\cr
  &\E k,\V i>k,\ \P a(x_i)\cr
  &\E k,\V i>k,\ \P(a,x_i)\cr
  &\PP(\aa,\<x_i>)\cr
  &\PP\aa(\<x_i>).\cr
}$$

\iind(viii) For $\<x_i>\in\UUn$, we have the equivalences:
$$\eqalign{
  &\EXT{(\V\D,\P)}(\<x_i>)\cr
  &\E k,\V i>k,\ (\V\D,\P)(x_i)\cr
  &\E k,\V i>k,\V y\in\D,\ \P(y,x_i)\cr
  &\V(y_i)\in\S\D,\E k,\V i>k,\ \P(y_i,x_i)\cr
  &\V\<y_i>\in\DD,\ \PP(\<y_i>,\<x_i>)\cr
  &(\V\DD,\PP)(\<x_i>).\cr
}$$
(It is easier to see that the fourth statement above implies the third by
negating both).

\iind(ix) If $\<x_i>\in\UUn$, then the assertions below are equivalent:
$$\eqalign{
  &\EXT{(\E\D,\P)}(\<x_i>)\cr
  &\E k,\V i>k,\ (\E\D,\P)(x_i)\cr
  &\E k,\V i>k,\E y\in\D,\ \P(y,x_i)\cr
  &\E(y_i)\in\S\D,\E k,\V i>k,\ \P(y_i,x_i)\cr
  &\E\<y_i>\in\DD,\ \PP(\<y_i>,\<x_i>)\cr
  &(\E\DD,\PP)(\<x_i>).\cr
}$$

\iind(x) Also, the following statements are equivalent:
$$\eqalign{
  &\EXT{(\EU\D,\P)}(\<x_i>)\cr
  &\E k,\V i>k,\ (\EU\D,\P)(x_i)\cr
  &\E k,\V i>k,\EU y\in\D,\ \P(y,x_i)\cr
  &\EU\<y_i>\in\DD,\ \PP(\<y_i>,\<x_i>)\cr
  &(\EU\DD,\PP)(\<x_i>).\cr
}$$

\iind(xi) As the following:
$$\eqalign{
  &\EXT{\P\o\f}(\<x_i>)\cr
  &\E k,\V i>k,\ (\P\o\f)(x_i)\cr
  &\E k,\V i>k,\ \P[\f(x_i)]\cr
  &\PP[\ff(\<x_i>)]\cr
  &(\PP\o\ff)(\<x_i>).\cr
}$$

\iind(xii) For any classes $\<x_i>\in\EXT{\D_f}$ and~$\<y_i>\in\EXT{\C_\g}$,
we have the equivalences:
$$\eqalign{
  &\EXT{\g\o\f}(\<x_i>)=\<y_i>\cr
  &\E k,\V i>k,\ (\g\o\f)(x_i)=y_i\cr
  &\E k,\V i>k,\ \g[\f(x_i)]=y_i\cr
  &\gg[\ff(\<x_i>)]=\<y_i>\cr
  &(\gg\o\ff)(\<x_i>)=\<y_i>.\cr
}$$
The equality $\EXT{\id\D}=\id\DD$ follows from item~(i) proved earlier.

\iind(xiii) If $\<x_i>\in\DD$ and
$\<y_i^j>=\<y_1^j,y_2^j,y_3^j,\ldots>\in\CC$ ($j=1,\ldots,n$), then these
statements are equivalent:
$$\eqalign{
  &\EXT{(\f_1,\ldots,\f_n)}(\<x_i>) = (\<y_i^1>,\ldots,\<y_i^n>)\cr
  &\E k,\V i>k,\ (\f_1,\ldots,\f_n)(x_i) = (y_i^1,\ldots,y_i^n)\cr
  &\E k,\V i>k,\ \f_j(x_i) = y_i^j \quad (j=1,\ldots,n)\cr
  &\EXT{f_j}(\<x_i>) = \<y_i^j> \quad (j=1,\ldots,n)\cr
  &(\EXT{\f_1},\ldots,\EXT{\f_n})(\<x_i>) = (\<y_i^1>,\ldots,\<y_i^n>),\cr
}$$
and also we have:
$$\eqalign{
  &\EXT{\p j\D}(\<y_i^1>,\ldots,\<y_i^n>) = \<x_i>\cr
  &\E k,\V i>k,\ \p j\D(y_i^1,\ldots,y_i^n) = x_i\cr
  &\E k,\V i>k,\ y_i^j = x_i\cr
  &\<y_i^j> = \<x_i>\cr
  &\p j\DD(\<y_i^1>,\ldots,\<y_i^n>) = \<x_i>.\cr
}$$

Finally, the last three items of the VET are corollaries of proposition~II.2:

\iind(a) By II.2(i), we have the equivalences:
$$\eqalign{
  &\V x\in\D,\ \P(x)\cr
  &\D\sc\P\cr
  &\DD\sc\PP\cr
  &\V\x\in\DD,\ \PP(\x).\cr
}$$

\iind(b) Using II.2(ii) and (iii) above, we see the following assertions are
equivalent:
$$\eqalign{
  &\E x\in\D,\ \P(x)\cr
  &\D\cap\P\ne\0\cr
  &\EXT{\D\cap\P}\ne\0\cr
  &\DD\cap\PP\ne\0\cr
  &\E\x\in\DD,\ \PP(\x).\cr
}$$

\iind(c) Now using II.2(iii) and item (iii), we have the equivalences:
$$\eqalign{
  &\EU x\in\D,\ \P(x)\cr
  &\D\cap\P \hbox{ is unitary}\cr
  &\EXT{\D\cap\P} \hbox{ is unitary}\cr
  &\DD\cap\PP \hbox{ is unitary}\cr
  &\EU\x\in\DD,\ \PP(\x).\cr
}$$

\tit{VI. Extension of Relation Attributes}

Our objective in this section is to illustrate the application of the~VET
with some basic examples. For that, let $\A$ be any subset of~$\Un$, where
$\U$ is the set from Sec.~IV and $n\in\N$ a natural number.

{\gr A binary relation $\P\sc\A^2$ is reflexive \IFF\ its virtual extension
$\PP\sc\AA^2$ is reflexive.}

\PF By the VET, the following statements are equivalent:
$$\eqalign{
  &\V x\in\A,\ \P(x,x)\cr
  &\V x\in\A,\ [\P\o(\id\A,\id\A)](x)\cr
  &\V\x\in\AA,\ \EXT{\P\o(\id\A,\id\A)}(\x)\cr
  &\V\x\in\AA,\ \PP\o(\EXT{\id\A},\EXT{\id\A})(\x)\cr
  &\V\x\in\AA,\ \PP\o(\id\AA,\id\AA)(\x)\cr
  &\V\x\in\AA,\ \PP(\x,\x),\cr
}$$
then $\P$ is reflexive \IFF\ $\PP$ is reflexive.\FP

{\gr A binary relation $\P\sc\A^2$ is symmetric \IFF\ its virtual extension
$\PP\sc\AA^2$ is symmetric.}

\PF Again by the VET, the following statements are equivalent:
$$\eqalign{
  &\V(x,y)\in\A^2,\ [\P(x,y)\IMP\P(y,x)]\cr
  &\V(x,y)\in\P,\ \P(y,x)\cr
  &\V(x,y)\in\P,\ [\P\o(\p 2{\A^2},\p 1{\A^2})](x,y)\cr
  &\V(\x,\y)\in\PP,\ \EXT{[\P\o(\p 2{\A^2},\p 1{\A^2})]}(\x,\y)\cr
  &\V(\x,\y)\in\PP,\ [\PP\o(\EXT{\p 2{\A^2}},\EXT{\p 1{\A^2}})](\x,\y)\cr
  &\V(\x,\y)\in\PP,\ [\PP\o(\p 2{\AA^2},\p 1{\AA^2})](\x,\y)\cr
  &\V(\x,\y)\in\PP,\ \PP(\y,\x)\cr
  &\V(\x,\y)\in\AA^2,\ [\PP(\x,\y)\IMP\PP(\y,\x)],\cr
}$$
so $\P$ is symmetric \IFF\ $\PP$ is symmetric.\FP

{\gr A binary relation $\P\sc\A^2$ is transitive \IFF\ its virtual extension
$\PP\sc\AA^2$ is transitive.}

\PF Once more the VET gives us:
$$\eqalign{
  &\V(x,y,z)\in\A^3,\ \{[\P(x,y)\AND\P(y,z)]\IMP\P(x,z)\}\cr
  &\V(x,y,z)\in[\P\o(\p 1{\A^3},\p 2{\A^3}) \AND
                \P\o(\p 2{\A^3},\p 3{\A^3})],\
               [\P\o(\p 1{\A^3},\p 3{\A^3})] (x,y,z)\cr
  &\V(\x,\y,\z)\in\EXT{\P\o(\p 1{\A^3},\p 2{\A^3}) \AND
                       \P\o(\p 2{\A^3},\p 3{\A^3})},\
                  \EXT{\P\o(\p 1{\A^3},\p 3{\A^3})} (\x,\y,\z)\cr
  &\V(\x,\y,\z)\in[\PP\o(\p 1{\AA^3},\p 2{\AA^3}) \AND
                   \PP\o(\p 2{\AA^3},\p 3{\AA^3})],
                   \ [\PP\o(\p 1{\AA^3},\p 3{\AA^3})] (\x,\y,\z)\cr
  &\V(\x,\y,\z)\in\AA^3,\ \{[\PP(\x,\y)\AND\PP(\y,\z)]\IMP\PP(\x,\z)]\},\cr
}$$
then $\P$ is transitive \IFF\ $\PP$ is transitive.\FP

Thus we have: {\gr a binary relation $\P\sc\A^2$ is an equivalence relation
on~$\A$ \IFF\ its virtual extension $\PP$ is an equivalence relation
on~$\AA$.}

We will say that a binary relation $\P\sc\A^2$ is {\gr antisymmetric\/} when:
$$
    \V(x,y)\in\A^2,\ [\P(x,y)\AND\P(y,x)]\IMP x=y.
$$
In addition, we will say that $\P$ is a {\gr partial ordering\/} on~$\A$ when
it is reflexive, transitive and antisymmetric. According to that, we have:

{\gr A binary relation $\P\sc\A^2$ is antisymmetric \IFF\ $\PP$ is
antisymmetric.}

\PF By the VET, the following statements are equivalent:
$$\eqalign{
  &\V(x,y)\in\A^2,\ [\P(x,y)\AND\P(y,x)]\IMP x=y\cr
  &\V(x,y)\in\{\P\AND[\P\o(\p 2{\A^2},\p 1{\A^2})]\},\ \eq\A(x,y)\cr
  &\V(\x,\y)\in\EXT{\P\AND[\P\o(\p 2{\A^2},\p 1{\A^2})]},\
               \EXT{\eq\A}(\x,\y)\cr
  &\V(\x,\y)\in\{\PP\AND[\PP\o(\p 2{\AA^2},\p 1{\AA^2})]\},\ \eq\AA(\x,\y)\cr
  &\V(\x,\y)\in\AA^2,\ [\PP(\x,\y)\AND\PP(\y,\x)]\IMP \x=\y,\cr
}$$
then $\P$ is antisymmetric \IFF\ $\PP$ is antisymmetric.\FP

So, we conclude that {\gr a binary relation $\P\sc\A^2$ is a partial ordering
on~$\A$ \IFF\ its virtual extension $\PP$ is a partial ordering on~$\AA$.}

A binary relation $\P\sc\A^2$ will be called a {\gr total ordering\/} when it
is reflexive, transitive, antisymmetric and satisfies the {\gr trichotomy:}
$$
    \V(x,y)\in\A^2,\ \P(x,y)\OR\P(y,x).
$$

Applying the VET to the trichotomy, as we did above, we obtain the following
equi\-val\-ences:
$$\eqalign{
  &\V(x,y)\in\A^2,\ [\P(x,y)\OR\P(y,x)](x,y)\cr
  &\V(x,y)\in\A^2,\ [\P\OR\P\o(\p2{\A^2},\p1{\A^2})](x,y)\cr
  &\V(\x,\y)\in\AA^2,\ \EXT{\P\OR\P\o(\p2{\A^2},\p1{\A^2})}(\x,\y).\cr
}$$
However, we cannot proceed as we did before because it is not true that the
extension of a logical disjunction is the same as the disjunction of the
extensions of disjunctives.

For instance, let us consider the total ordering `$\le$' between real
numbers.  The extension `$\EXT\le$' is just a partial ordering relation
between virtual numbers. If $\a\in\RR$ is the class of the sequence
$(-1,+1,-1,+1,\ldots)\in\S\R$, then the statements `$\a\EXT\le\EXT0$'
and~`$\EXT0\EXT\le\a$' are both false.

In the case of the connectives `not', `$\IMP$', and~`$\EQV$' we have a
similar situation, but in the opposite direction, since inclusions present in
items (ii), (v) and~(vi) of the~VET are opposite the one occurring in
item~(iv). In spite of that, the VET can establish unidirectional
implications between statements involving those connectives, as shown in the
following examples:

(i) If the extension $\PP$ satisfies the trichotomy then $\P$ also
satisfies it, since, by item~(iv) of the~VET, the condition:
$$
    \EXT{\P\OR\P\o(\p 2{\A^2},\p 1{\A^2})}(\x,\y)
$$
is necessary (although not sufficient) for the validity of:
$$
    [\PP\OR\EXT{\P\o(\p 2{\A^2},\p 1{\A^2})}](\x,\y).
$$

(ii) If $\P$ and~$\Q$ are two binary relations such that:
$$
    \V y\in\A,\E x\in\A,\ \P(x,y)\IMP\Q(x,y),
$$
then:
$$
    \V\y\in\AA,\E\x\in\AA,\ \PP(\x,\y)\IMP\QQ(\x,\y)
$$
holds. To verify this, we initially have the equivalences:
$$\eqalign{
  &\V y\in\A,\E x\in\A,\ \P(x,y)\IMP\Q(x,y)\cr
  &\V y\in\A,\E x\in\A,\ (\P\IMP\Q)(x,y)\cr
  &\V y\in\A,\ [\E\A,(\P\IMP\Q)](y)\cr
  &\V\y\in\AA,\ \EXT{\E\A,(\P\IMP\Q)}(\y)\cr
  &\V\y\in\AA,\ [\E\AA,\EXT{\P\IMP\Q}](\y)\cr
  &\V\y\in\AA,\E\x\in\AA,\ \EXT{\P\IMP\Q}(\x,\y).\cr
}$$
Furthermore, by item~(v) of the~VET, the condition:
$$
    \EXT{\P\IMP\Q}(\x,\y)
$$
is sufficient to guarantee that:
$$
    (\PP\IMP\QQ)(\x,\y).
$$
Therefore, the last statement of the series above implies
$$
    \V\y\in\AA,\E\x\in\AA,\ (\PP\IMP\QQ)(\x,\y),
$$
which is equivalent to:
$$
    \V\y\in\AA,\E\x\in\AA,\ \PP(\x,\y)\IMP\QQ(\x,\y).
$$

\tit{VII. Extension of Function and Operation Properties}

In this section we will introduce some basic applications of the~VET
involving functions and operations. Let $\D\sc\Un$ and $\C\sc\Uk$ be two
generic subsets, where $\U$ is the set considered in Sec.~IV, and $n,k\in\N$
two natural numbers.

{\gr A function $\f\:\D\to\C$ is one-to-one \IFF\ its virtual extension
$\ff\:\DD\to\CC$ is one-to-one.}

\PF By the VET, the following statements are equivalent:
$$\eqalign{
  &\V(x_1,x_2)\in\D^2,\ [\f(x_1)=\f(x_2)\ \IMP\ x_1=x_2]\cr
  &\V(x_1,x_2)\in[\eq\C\o(\f\o\p 1{\D^2},
     \f\o\p 2{\D^2})],\ \eq\D(x_1,x_2)\cr
  &\V(\x_1,\x_2)\in\EXT{\eq\C\o(\f\o\p 1{\D^2},
     \f\o\p 2{\D^2})},\ \EXT{\eq\D}(\x_1,\x_2)\cr
  &\V(\x_1,\x_2)\in[\eq\CC\o(\ff\o\p 1{\DD^2},
     \f\o\p 2{\DD^2})],\ \eq\DD(\x_1,\x_2)\cr
  &\V(\x_1,\x_2)\in\DD^2,\ [\f(\x_1)=\f(\x_2)\ \IMP\ \x_1=\x_2],\cr
}$$
then $\f$ is one-to-one \IFF\ $\ff$ is one-to-one.\FP

{\gr A function $\f$ maps $\D$ onto $\C$ \IFF\ its virtual extension $\ff$
maps $\DD$ onto~$\CC$.}

\vbox{\PF Again by the VET, the following statements are equivalent:
\medskip
$$\eqalign{
  &\V y\in\C,\E x\in\D,\ \f(x)=y\cr
  &\V y\in\C,\E x\in\D,\ \eq\C\o(\f\o\p 1{\D\X\C},\p 2{\D\X\C})(x,y)\cr
  &\V y\in\C,\ [\E\D,\eq\C\o(\f\o\p 1{\D\X\C},\p 2{\D\X\C})](y)\cr
  &\V\y\in\CC,\ \EXT{\E\D,\eq\C\o(\f\o\p 1{\D\X\C},\p 2{\D\X\C})}(\y)\cr
  &\V\y\in\CC,\ [\E\DD,\eq\CC\o(\ff\o\p 1{\DD\X\CC},\p 2{\DD\X\CC})](\y)\cr
  &\V\y\in\CC,\E\x\in\DD,\ \eq\CC\o(\ff\o\p 1{\DD\X\CC},\p 2{\DD\X\CC})
      (\x,\y)\cr
  &\V\y\in\CC,\E\x\in\DD,\ \ff(\x)=\y.\FP\cr
}$$}
\medskip

Therefore, {\gr a function $\f\:\D\to\C$ is inversible \IFF\ its virtual
extension $\ff\:\DD\to\CC$ is inversible}. In this case, it follows directly
from item~(xii) of the~VET that the extension of the inverse function of~$\f$
is equal to the inverse function of its virtual extension.

Let now $\op\:\A^2\to\A$ be a binary operation on~$\A\sc\Un$. Applying the
VET as we did above, we verify that {\gr $\op$ is associative \IFF\
$\EXT{\op}$ is associative}, and that {\gr $\op$ is commutative \IFF\
$\EXT{\op}$ is commutative.} Besides, if $\od\:\A^2\to\A$ is another binary
operation defined on the same set $\A$, then {\gr $\od$ is distributive with
respect to~$\op$ \IFF\ $\EXT\od$ is distributive with respect to~$\EXT\op$}.
In other words, the equality:
\medskip
$$
    a\od(b\op c) = (a\od b)\op(a\od c)
$$
\medskip\noindent
holds for every triple $(a,b,c)\in\A^3$ \IFF\
\medskip
$$
    \a\EXT\od(\b\EXT\op\ga) = (\a\EXT\od\b)\EXT\op(\a\EXT\od\ga)
$$
\medskip\noindent
holds for every triple $(\a,\b,\ga)\in\AA^3$.

As to the neutral element, we have:

{\gr If $\op\:\A^2\to\A$ is a binary operation then $e\in\A$ is a right
neutral element of~$\op$ \IFF\ $\EXT e\in\KA\sc\AA$ is a right neutral
element of~$\EXT\op$}.

\vbox{\PF By the VET, the following statements are equivalent:
\smallskip
$$\eqalign{
  &\V a\in\A,\ e\op a = a\cr
  &\V a\in\A,\ [\eq\A\o(\op,\p2{\A^2})](e,a)\cr
  &\V a\in\A,\ [\eq\A\o(\op,\p2{\A^2})]e(a)\cr
  &\V\a\in\AA,\ \EXT{[\eq\A\o(\op,\p2{\A^2})]e}(\a)\cr
  &\V\a\in\AA,\ [\eq\AA\o(\EXT\op,\p2{\AA^2})]\EXT e(\a)\cr
  &\V\a\in\AA,\ [\eq\AA\o(\EXT\op,\p2{\AA^2})](\EXT e,\a)\cr
  &\V\a\in\AA,\ \EXT e\EXT\op\a = \a.\FP\cr
}$$}
\smallskip

Analogously:

{\gr If $\op\:\A^2\to\A$ is a binary operation then $e\in\A$ is a left
neutral element of~$\op$ \IFF\ $\EXT e\in\KA\sc\AA$ is a left neutral element
of~$\EXT\op$}.

\PF Once more the VET gives us:
\smallskip
$$\eqalign{
  &\V a\in\A,\ a\op e = a\cr
  &\V a\in\A,\ \{\eq\A\o[\op\o(\p2{\A^2},\p1{\A^2}),\p2{\A^2}]\}(e,a)\cr
  &\V a\in\A,\ \{\eq\A\o[\op\o(\p2{\A^2},\p1{\A^2}),\p2{\A^2}]\}e(a)\cr
  &\V\a\in\AA,\ \EXT{\{\eq\A\o
      [\op\o(\p2{\A^2},\p1{\A^2}),\p2{\A^2}]\}e}(\a)\cr
  &\V\a\in\AA,\ \{\eq\AA\o
      [\EXT\op\o(\p2{\AA^2},\p1{\AA^2}),\p2{\AA^2}]\}\EXT e(\a)\cr
  &\V\a\in\AA,\ \{\eq\AA\o
      [\EXT\op\o(\p2{\AA^2},\p1{\AA^2}),\p2{\AA^2}]\}(\EXT e,\a)\cr
  &\V\a\in\AA,\ \a\op\EXT e = \a.\FP\cr
}$$
\bigskip

Now we will consider the existence of opposites:

{\gr If $\op\:\A^2\to\A$ and~$c\in\A$, then the condition:
$$
    \V a\in\A,\E b\in\A,\ a\op b = c
$$
is necessary and sufficient to the validity of:}
$$
    \V\a\in\AA,\E\b\in\AA,\ \a\EXT\op\b = \EXT c.
$$

\PF By the VET, the following statements are equivalent:
$$\eqalign{
  &\V a\in\A,\E b\in\A,\ a\op b = c\cr
  &\V a\in\A,\E b\in\A,\ \{\eq\A\o[\op\o(\p2{\A^3},\p3{\A^3}),\p1{\A^3}]\}
      (c,a,b)\cr
  &\V a\in\A,\E b\in\A,\ \{\eq\A\o[\op\o(\p2{\A^3},\p3{\A^3}),\p1{\A^3}]\}
      c(a,b)\cr
  &\V\a\in\AA,\E\b\in\AA,\ \EXT{\{\eq\A\o
      [\op\o(\p2{\A^3},\p3{\A^3}),\p1{\A^3}]\}c}(a,b)\cr
  &\V\a\in\AA,\E\b\in\AA,\ \{\eq\AA\o
      [\EXT\op\o(\p2{\AA^3},\p3{\AA^3}),\p1{\AA^3}]\}\EXT c(\a,\b)\cr
  &\V\a\in\AA,\E\b\in\AA,\ \{\eq\AA\o
      [\EXT\op\o(\p2{\AA^3},\p3{\AA^3}),\p1{\AA^3}]\}(\EXT c,\a,\b)\cr
  &\V\a\in\AA,\E\b\in\AA,\ \a\EXT\op\b = \EXT c.\FP\cr
}$$

Nevertheless, the condition:
$$
    \V a\neq d,\E b\in\A,\ a\od b = c
$$
is equivalent to:
$$
    \V\a\EXT\neq\EXT d,\E\b\in\AA,\ \a\EXT\od\b = \EXT c,
$$
which, by item~(ii) of the~VET, is necessary but not sufficient to assure the
validity of:
$$
    \V\a\neq\EXT d,\E\b\in\AA,\ \a\EXT\od\b = \EXT c.
$$

\tit{VIII. Extension of Mathematical Structures}

The results of the last two sections illustrate how the VET can be used to
logically relate a statement about the set~$\U$ to another statement about
its virtual extension. The following syntactic rules informally describe how
that ``extended statement'' is obtained from the original:

(i) consistently substitute the bound variables (quantified) ranging over a
subset $\A\sc\Un$ by bound variables ranging over its virtual extension
$\AA\sc\UUn$, keeping the corresponding quantifier;

(ii) consistently substitute the free variables (not quantified) by the
corresponding element in $\K\U\sc\UU$;

(iii) substitute the functions present in the original statement by the
respective virtual extensions;

(iv) selectively substitute the relations in the original statement by its
virtual extensions, respecting the restrictions on the connectives `not',
`or', `$\IMP$' and~`$\EQV$'.

However, it is important to note that the~VET has been enunciated and proved
by methods of Elementary Set Theory, not having used the formal distinction
between syntactic and semantic planes which characterizes Mathematical Logic.

The results presented in the last two sections show that application of
the~VET also does not require more than elementary mathematical techniques.
Collecting some of those results, we conclude that:

{\gr A pair $(\G,\op)$ is a group \IFF\ its virtual extension $(\EXT
G,\EXT\op)$ is a group. In this case, $(\G,\op)$ is commutative \IFF\ $(\EXT
G,\EXT\op)$ is commutative.}

{\gr A triple $(\A,\op,\od)$ is a ring \IFF\ its virtual extension $(\EXT
G,\EXT\op,\EXT\od)$ is a ring. In this case, $e\in\A$ is a unity in~$\A$
\IFF\ $\EXT e\in\KA\sc\AA$ is a unity in~$\AA$.}

Nevertheless, virtual extension of a total ordering is just partial ordering,
and the virtual extension of a field is just a ring with unity. For example,
the virtual extension $(\RR,\EXT+,\EXT\X,\EXT<)$ of the ordered field of real
numbers is not a field, nor is it totally ordered.

On the other hand, every relation and function defined on any set $\A$ is
extended to the set~$\AA$. In addition, the loss of part of the mathematical
structure of~$\A$ during the process of virtual extension can be compensated
by the~VET, which allows us to transport many facts about those extended
relations and functions directly to~$\AA$. For instance, we have the
trigonometric functions:
$$
    \EXT\sin\:\RR\to\RR \qquad\hbox{and}\qquad \EXT\cos\:\RR\to\RR,
$$
which satisfy the identity
$$
    \EXT\sin^2\a\EXT+\EXT\cos^2\a = \EXT1,
$$
for any virtual number $\a\in\RR$.

An important aspect of the virtual extension process is that it must be
applied ``sim\-ul\-taneously'' to every set bound by the relations which we
intend to extend. We can do that just by taking the ``disjoint union'' of
these sets as ``universe''. In other words, if we intend to extend relations
between the sets of a family $\A_i$ ($i\in\I$), then we should take the
disjoint union of the family $(\A_i)$ as the set represented by~$\U$ in the
previous sections.

As an illustration, to extend a vectorial space $\Ve$ over real scalars, we
can make $\U$ equal the disjoint union of $\Ve$ and~$\R$, so that a generic
member of~$\UU$ will be a class of sequences whose elements can be either
vectors of~$\Ve$ or real scalars. In this manner, the {\gr virtual vectors\/}
will be the members of the subset $\VV\sc\UU$ (i.e., the classes of sequences
which end taking values only in~$\Ve$), whereas the {\gr virtual scalars\/}
will be the members of the subset $\RR\sc\UU$ (which is the class of
sequences that end taking only real values). Thus, the multiplication of
scalars and vectors:
$$
    \X\:\R\X\Ve \to \Ve
$$
extends to an operation between virtual scalars and virtual vectors:
$$
    \EXT\X\:\RR\X\VV \to \VV.
$$

Another example of mathematical structure whose definition involves more than
one set is that of {\gr manifold\/} (topological or differentiable).
Proceeding as above, we can construct a {\gr virtual manifold\/} armed with
an atlas of {\gr virtual charts}, which associate {\gr local virtual
coordinates\/} to the {\gr virtual points\/} of the manifold.

The word ``extension'' is commonly used in mathematics to indicate that we
have a ``copy'' of the original set $\U$ inside the set $\UU$ constructed
from it, and that we intend to identify $\U$ with that copy. The~VET
guarantees that $\K\U\sc\UU$ is a faithful copy of~$\U$, since any relation
$\P$ involving the members of~$\U$ is equivalent to the restriction of its
extension $\PP$ to the corresponding members in~$\K\U$. This fact authorizes
the identification $\U=\K\U$, which allows us to consider $\U\sc\UU$. If $\U$
is the disjoint union of family $(\A_i)$, then we will have
$\A_i=\K{\A_i}\sc\EXT{\A_i}\sc\UU$, for every $i\in\I$.
\bye